\documentclass[12pt]{amsart}

\usepackage{upref}

\newcommand{\integer}{{\mathord{\mathbb{Z}}}}
\newcommand{\rational}{{\mathord{\mathbb{Q}}}}
\newcommand{\real}{{\mathord{\mathbb{R}}}}
\newcommand{\complex}{{\mathord{\mathbb{C}}}}
\newcommand{\quaternion}{{\mathord{\mathbb{H}}}}
\newcommand{\Qi}{{\rational(i)}}

\newcommand{\otimesR}{\mathbin{\otimes_\real}}
\newcommand{\otimesQ}{\mathbin{\otimes_\rational}}
\newcommand{\otimesF}{\mathbin{\otimes_F}}

\newcommand{\Qrank}{\mathop{\mbox{\upshape$\rational$-rank}}}
\newcommand{\Rrank}{\mathop{\mbox{\upshape$\real$-rank}}}

\newcommand{\Lie}[1]{\mathfrak{\lowercase{#1}}}
\newcommand{\conj}[1]{\overline{#1}}
\newcommand{\iso}{\cong}
\newcommand{\Mat}{\operatorname{Mat}\nolimits}
\newcommand{\End}{\operatorname{End}\nolimits}

\newcommand{\Id}{\mathord{\mathrm{Id}}}
\newcommand{\SU}{\operatorname{SU}}
\newcommand{\su}{\mathop{\Lie{SU}}}
\newcommand{\so}{\operatorname{\Lie{so}}}
\newcommand{\symp}{\operatorname{\Lie{sp}}}

\newcommand{\closure}[1]{\overline{#1}}

\newcommand{\cjg}[1]{\overline{#1}}

\renewcommand{\see}[1]{{\upshape(}see~\ref{#1}{\upshape)}}
\newcommand{\cf}[1]{{\upshape(}cf.~\ref{#1}{\upshape)}}
\newcommand{\fullsee}[2]{{\upshape(}see~\ref{#1}\pref{#1-#2}{\upshape)}}

\newcommand{\seeSect}[1]{{\upshape(}see~\S\ref{#1}{\upshape)}}

\newcommand{\pref}[1]{{\upshape(}\ref{#1}{\upshape)}}
\newcommand{\fullref}[2]{{\ref{#1}\pref{#1-#2}}}

\swapnumbers

\newtheorem{thm}[equation]{Theorem}
\newtheorem{prop}[equation]{Proposition}
\newtheorem{cor}[equation]{Corollary}
\newtheorem{lem}[equation]{Lemma}

\theoremstyle{definition}
\newtheorem{defn}[equation]{Definition}
\newtheorem{notation}[equation]{Notation}
\newtheorem{rem}[equation]{Remark}
\newtheorem{eg}[equation]{Example}
\newtheorem{ack}[equation]{Acknowledgments}

\numberwithin{equation}{section}

 \newcounter{step}

\newenvironment{step}[1][\unskip]{\refstepcounter{step}
 \setcounter{substep}{0}
 \em
 \medskip \noindent Step \thestep\ #1.\ }{\unskip\upshape}
 \renewcommand{\thestep}{\arabic{step}}

 \newcounter{substep}

\newenvironment{substep}[1][\unskip]{\refstepcounter{substep}\em
 \medskip \noindent Substep \thesubstep\ #1.\
}{\unskip\upshape}
 \renewcommand{\thesubstep}{\thestep.\arabic{substep}}

 \newcounter{case}

\newenvironment{case}[1][\unskip]{\refstepcounter{case}
 \setcounter{subcase}{0}\em
 \medskip \noindent Case \thecase\ #1.\ }{\unskip\upshape}
 \renewcommand{\thecase}{\arabic{case}}

 \newcounter{subcase}
 \newenvironment{subcase}[1][\unskip]{\refstepcounter{subcase}\em
 \medskip \noindent Subcase \thesubcase\ #1.\ }{\unskip\upshape}
 \renewcommand{\thesubcase}{\arabic{case}.\arabic{subcase}}

\begin{document}

\title[Real representations have $\rational$-forms]
 {Real representations of semisimple Lie algebras
have $\rational$-forms}

\author{Dave Witte}

\address{Department of Mathematics,
 Oklahoma State University,
 Stillwater, OK 74078, USA}

\email{dwitte@math.okstate.edu, 
  \hfil\break
 http://www.math.okstate.edu/$\sim$dwitte}

\dedicatory{To Professor M.~S.~Raghunathan on his sixtieth
birthday}

\date{17 December 2002} 

\begin{abstract}
 We prove that each real semisimple Lie algebra~$\Lie G$ has a
$\rational$-form $\Lie G_\rational$, such that every real
representation of~$\Lie G_\rational$ can be realized
over~$\rational$. This was previously proved by
M.~S.~Raghunathan (and rediscovered by P.~Eberlein) in the
special case where $\Lie G$ is compact.
 \end{abstract}

\thanks{Submitted in May 2002 to the proceedings of the
conference
 \emph{Algebraic Groups and Arithmetic {\upshape(}Mumbai, India,
17-22 December 2001\upshape{)}}. Revised December 2002.} 

\maketitle

\section{Introduction}

All Lie algebras and all representations are assumed to be
finite-dimensional. It is easy to see, from the theory of
highest weights, that if $\Lie G$ is an $\real$-split,
semisimple Lie algebra over~$\real$, then every
$\complex$-representation of~$\Lie G$ has an $\real$-form
\see{split->reps}. (That is, if $V_\complex$ is a
representation of~$\Lie G$ over~$\complex$, then there is a
real representation~$V$ of~$\Lie G$, such that $V_\complex
\iso V \otimesR \complex$.) Because every semisimple
Lie algebra over~$\complex$ has an $\real$-split real form,
this leads to the following immediate conclusion:

\begin{rem}
   Any complex semisimple Lie algebra~$\Lie G_\complex$ has a
real form~$\Lie G$, such that every $\complex$-representation
of~$\Lie G$ has a real form.
 \end{rem}

 In this paper, we prove the analogous statement with the
field extension $\complex/\real$ replaced with
$\real/\rational$.

\begin{thm}[{see \fullref{EGoodQform}{Qform}}]
  Any real semisimple Lie algebra~$\Lie G$ has a
$\rational$-form $\Lie G_\rational$, such that every real
representation of~$\Lie G_\rational$ has a $\rational$-form.
 \end{thm}

 In the special case where $\Lie G$ is compact, the theorem was
proved by M.~S.~Raghunathan \cite[\S3]{Raghunathan-Qform}. This
special case was independently rediscovered by P.~Eberlein
\cite{EberleinQform}, and a very nice proof was found by R.~Pink
and G.~Prasad (personal communication,
see~\S\ref{PinkPrasadMethSect}). When $\Lie G$ is compact, these
authors showed that the ``obvious" $\rational$-form of~$\Lie G$
has the desired property.

At the other extreme, where $\Lie G$ is $\real$-split, we may
take $\Lie G_\rational$ to be any $\rational$-split
$\rational$-form of~$\Lie G$ \see{split->reps}. 

 The general case is a combination of the two extremes, and
the desired $\rational$-form can be obtained from a Chevalley
basis of $\Lie G \otimesR \complex$ by slightly modifying a
construction of A.~Borel \cite{Borel-CK}
\seeSect{ConstructQformSect}. We give two different proofs that
this $\rational$-form has the desired property: one proof is by
the method of Pink and Prasad, using a little bit of number
theory \seeSect{PinkPrasadMethSect}, and the other proof is by
reducing to the compact case, so Raghunathan's theorem applies
\seeSect{RaghuMethodSect}.

It would be interesting to characterize the semisimple Lie
algebras $\Lie G_\rational$ over~$\rational$, such that every
real representation has a $\rational$-form. For example, work of
J.~Tits \cite{Tits-IrredRepns} implies that every
$\rational$-form of $\symp(n)$ has this property
\see{TitsAllButSO}. On the other hand, it is important to note
that there exist examples of $\Qi$-split Lie algebras that do
not have this property \see{badGQ}. (Real
representations of such a Lie algebra can be realized over both
$\Qi$ and~$\real$, but not over $\Qi \cap \real  = \rational$.)

\begin{ack}
 I am very grateful to Gopal Prasad for sharing with me the
elegant proof that he found in collaboration with R.~Pink, and
for his many other helpful comments on the original version of
this manuscript.

 I would like to thank Bob Stanton, for bringing the work of
P.~Eberlein to my attention, Scot
Adams and Patrick Eberlein, for encouraging me to think more
carefully about this problem and for their many helpful
comments, T.~N.~Venkataramana, for sharing his
insights into the material of~\S\ref{ConstructQformSect}, and
Nilabh Sanat, for explaining part of Lem.~\ref{RealVsQuat} to me.
 The research was partially supported by grants from the
National Science Foundation (DMS--9801136 and DMS--0100438) and
the German-Israeli Foundation for Research and Development. 

 Much of the work was carried out during visits to the Duluth
campus of the University of Minnesota and to the University of
Bielefeld. The manuscript was revised into a publishable form
at the University of Chicago, and some later revisions were made
at the University of Lethbridge. I am pleased to thank the
mathematics departments of all four of these institutions for
their hospitality that made these visits so productive, and so
enjoyable.
 \end{ack}

\section{More precise statement of the main result}

\begin{prop}[Pink, Prasad, see \S\ref{PinkPrasadMethSect}]
\label{Quasi->good}
 Suppose $G$ is a connected, reductive algebraic
$\rational$-group. If
 \begin{enumerate} \renewcommand{\theenumi}{\alph{enumi}}
 \item $G$ is split over some imaginary quadratic
extension~$F$ of~$\rational$,
 and
 \item \label{Quasi->good-quasi}
 $G$ is quasi-split over the $p$-adic field $\rational_p$,
for every odd prime~$p$,
 \end{enumerate}  
 then each irreducible $\rational$-representation of~$G$
remains irreducible over~$\real$.
 \end{prop}

This can be restated in the following equivalent form.

\begin{defn}
 Suppose $(\pi, V)$ is a real representation of an algebraic
$\rational$-group~$G$. 
 A $\rational$-subspace $V_\rational$ of~$V_\real$ is a
\emph{$\rational$-form} of $(\pi, V)$ if 
 \begin{itemize}
 \item $V_\rational$ is $G_\rational$-invariant,
 and
 \item $V_\rational$ is the $\rational$-span of an
$\real$-basis of~$V_\real$ (so $V_\real \iso V_\rational
\otimesQ \real$).
 \end{itemize}
 \end{defn}

\begin{cor}[see~\ref{Qform<>RIrred}] \label{RrepsQform}
 If $G$ is as in Prop.~\ref{Quasi->good}, then every real
representation of~$G$ has a $\rational$-form.
 \end{cor}

\begin{prop}[see \S\ref{ConstructQformSect}]
\label{GRhasQform}
 Every connected, simply connected, semisimple real algebraic
group has a $\rational$-form satisfying the hypotheses of
Prop.~\ref{Quasi->good}.
 \end{prop}

Combining Prop.~\ref{GRhasQform} with
Prop.~\ref{Quasi->good} and Cor.~\ref{RrepsQform} immediately
yields the following conclusion.

\begin{defn}
 Suppose
 \begin{itemize}
 \item $\Lie G_\rational$ is a Lie algebra over~$\rational$,
 and
 \item $(\pi, V)$ is a real representation of~$\Lie
G_\rational$.
 \end{itemize}
 A $\rational$-subspace $V_\rational$ of~$V$ is a
\emph{$\rational$-form} of $(\pi, V)$ if 
 \begin{itemize}
 \item $V_\rational$ is $\Lie G_\rational$-invariant,
 and
 \item $V_\rational$ is the $\rational$-span of an
$\real$-basis of~$V$ (so $V \iso V_\rational
\otimesQ \real$).
 \end{itemize}
 \end{defn}

\begin{cor} \label{EGoodQform}
 Any real semisimple Lie algebra~$\Lie G$ has a
$\rational$-form $\Lie G_\rational$, such that
 \begin{enumerate}
 \item \label{EGoodQform-irred}
 if $V_\rational$ is any irreducible
$\rational$-representation of~$\Lie G_\rational$, then the
$\real$-representation $V_\real = V_\rational
\otimesQ \real$ is irreducible; and
 \item \label{EGoodQform-Qform}
 every real representation of~$\Lie G_\rational$ has a
$\rational$-form.
 \end{enumerate}
 \end{cor}

See Thm.~\ref{Weyl->good} for a version of
Prop.~\ref{Quasi->good} that replaces
\fullref{Quasi->good}{quasi} with a quite different hypothesis,
due to M.~S.~Raghunathan. The $\rational$-form constructed in
\S\ref{ConstructQformSect} also satisfies this alternate
hypothesis, so this yields a different proof of
Cor.~\ref{EGoodQform}.

\section{Preliminaries}

The following is well known (see, for example,
\cite[Thm.~2.5]{Tits-IrredRepns}).

\begin{lem} \label{split->reps}
 Let
 \begin{itemize}
 \item $F$ be a subfield of~$\complex$,
 \item $\Lie G$ be a semisimple Lie algebra over~$F$,
 and
 \item $V_\complex$~be a $\complex$-representation
of~$\Lie G$.
 \end{itemize}
 If $\Lie G$ is $F$-split, then $V$ has an $F$-form.
 \end{lem}

\begin{proof}
 Let $\Lie T$ be a maximal $F$-split torus of~$\Lie G$. Because
every representation of~$\Lie G$ is a direct sum of
irreducibles, we may assume $V_\complex$ is irreducible; let
$\lambda$~be the highest weight of~$V_\complex$ (with respect to
some ordering of the roots of~$\Lie T$). Since $\lambda$ is a
character of the $F$-split torus~$\Lie T$, we know that
$\lambda(\Lie T_F) \subset F$. So there is an
$F$-representation~$V_F$ of~$\Lie G$ with highest
weight~$\lambda$. Hence, $V_F \mathbin{\otimes_F} \complex \iso
V_\complex$, so $V_F$ is (isomorphic to) an $F$-form
of~$V_\complex$.
 \end{proof}

For future reference, let us record the following consequence of
this fact.

\begin{cor} \label{Q(i)split}
 Suppose, for some quadratic extension~$F$ of~$\rational$, that
$\Lie G$ is an $F$-split, semisimple Lie algebra
over~$\rational$.
 \begin{enumerate}
 \item \label{Q(i)split-irred}
 If $V_F$ is any irreducible $F$-representation of~$\Lie G_F$,
then $V_\complex$ is irreducible.
 \item \label{Q(i)split-2irreds}
 If $V_\rational$~is any irreducible $\rational$-representation
of~$\Lie G$, then $V_\complex$ is either irreducible or the
direct sum of two irreducibles.
 \end{enumerate}
 \end{cor}

\begin{proof}
 \pref{Q(i)split-irred} The proof of
Lem.~\ref{Qform<>RIrred}($\Leftarrow$) below shows that this is a
consequence of Lem.~\ref{split->reps}.

\pref{Q(i)split-2irreds} Write $F = \rational[\sqrt{r}]$, for
some $r \in \rational$. Because $V_F = V_\rational + \sqrt{r}
V_\rational$ (and $V_\rational$ is an irreducible $\Lie
G_\rational$-module), we know that the $\Lie G_F$-module $V_F$
is either irreducible or the direct sum of two irreducibles.
Then the desired conclusion follows from~\pref{Q(i)split-irred}.
 \end{proof}

The following observation must be well known. The
direction ($\Rightarrow$) can be found in
\cite[\S3]{Raghunathan-Qform}.

\begin{lem} \label{Qform<>RIrred}
 Suppose $G$ is a connected, semisimple algebraic group
over~$\rational$. Every irreducible $\rational$-representation
of~$G$ remains irreducible over~$\real$ if and only if every
real representation of~$G$ has a $\rational$-form.
 \end{lem}

\begin{proof} 
 ($\Leftarrow$) Let $V$ be a $\rational$-representation
of~$G$, such that $V_\real$~is reducible. (We wish to show
that $V$~is reducible.) We may write $V_\real = U_1 \oplus
U_2$, for some nontrivial $\real$-representations $U_1$
and~$U_2$. By assumption, there exist
$\rational$-representations $V_1$ and~$V_2$, such that
$(V_j)_\real \iso U_j$. Then 
 $$V_\real = U_1 \oplus U_2 \iso (V_1)_\real \oplus
(V_2)_\real \iso (V_1 \oplus V_2)_\real .$$
 Thus, $V$ is isomorphic to $V_1 \oplus V_2$ over~$\real$.
Since both $V$ and $V_1 \oplus V_2$ are defined
over~$\rational$, this implies that $V$ is isomorphic to $V_1
\oplus V_2$ over~$\rational$. (The $\Lie G$-equivariant maps
from $V_\real$ to $(V_1 \oplus V_2)_\real$ form a real vector
space that is defined over~$\rational$, so the
$\rational$-points span.) Thus, $V$ is reducible.

($\Rightarrow$) Let $V$ be a real representation of~$G$. To
simplify notation (and because this is the only case we
need), let us assume that $G$~is split over some imaginary
quadratic extension~$F$ of~$\rational$. Then $V_\complex$ has an
$F$-form~$U$. Let $U|_\rational$ be the
$\rational$-representation obtained by viewing $U$ as a vector
space over~$\rational$. 

Now write $U|_\rational = U_1 \oplus \cdots \oplus U_r$ as a
direct sum of irreducible $\rational$-modules. Then 
 $$ V_\complex|_\real
 = U_\complex|_\real
 \iso (U|_\rational)_\real
 \iso (U_1)_\real \oplus \cdots \oplus (U_r)_\real $$
 
Since $V$ is a submodule of $V_\complex|_\real$ (indeed,
$V_\complex|_\real$ is the direct sum $V \oplus i V$ of two
copies of~$V$), and, by assumption, each $(U_j)_\real$ is
irreducible, we conclude that $V$ is isomorphic to
$(U_j)_\real$, for some~$j$. So (up to isomorphism) $U_j$~is a
$\rational$-form of~$V$.
 \end{proof}

We also use the following (special case of a) result of J.~Tits
\cite[Thms.~7.2(i) and 3.3]{Tits-IrredRepns} that applies to the
quasi-split case. In our applications, $F$~will be a $p$-adic
field~$\rational_p$ \fullsee{Quasi->good}{quasi}.

\begin{prop}[Tits] \label{QuasiSplit->Field}
 Let
 \begin{itemize}
 \item $F$ be a field of chacteristic zero,
 \item $G$ be a connected, reductive algebraic $F$-group,
 and
 \item $(\pi, V)$~be an irreducible $F$-representation
of~$G$.
 \end{itemize}
 If $G$ is quasi-split over~$F$, then $\End_G(V)$ is
commutative.
 \end{prop}

\begin{proof}
 By assumption, $G$ has a Borel subgroup~$B$ that is defined
over~$F$. Let $v \in V$ be a nonzero vector that is fixed by
every element of the unipotent radical of~$B$, and let
$\closure{F}$ be the algebraic closure of~$F$. 

Schur's Lemma asserts that $\End_G(V)$ is a division algebra.
By enlarging~$F$, we may assume that $F$ is the center of
$\End_G(V)$, so $\End_G(V_{\overline{F}}) = \End_G(V) \otimesF
\overline{F}$ is a simple algebra, which implies that
$V_{\overline{F}}$ is isotypic: we have $V_{\closure{F}} = W
\oplus \cdots \oplus W$, for some irreducible
$\Lie G_{\closure{F}}$-module~$W$. This implies that $v$~is a
weight vector (that is, an eigenvector for $\pi(B)$), so $V$ is a
highest-weight module. Therefore $\End_G(V) = F$, which
is abelian, as desired.
 \end{proof}

The method of R.~Pink and G.~Prasad utilizes the following
classical result of number theory. It is obtained by combining
the Hasse Principal (or local-to-global principle) with the fact
that the sum of the local invariants of a quaternion algebra (or,
what is the same thing, of a quadratic form) is~$0$ (so, if all
but one of them vanish, then they all must vanish). 

\begin{lem}[{cf.~\cite[Cor.~3.2.3, p.~43]{Serre-CourseArith}}]
\label{SplitLocal->Split}
 Let $\mathcal{C}$ be a quaternion division algebra
over~$\rational$. If $\mathcal{C}$ splits over~$\rational_p$,
for every odd prime~$p$, then $\mathcal{C}$ does \textbf{not}
split over~$\real$.
 \end{lem}

\begin{rem} \label{HasseAnyPrime}
 The exceptional prime~$2$ can be replaced with any other prime
in Lem.~\ref{SplitLocal->Split}: if there is a prime~$p_0$, such
that $G$ is quasi-split over the $p$-adic field $\rational_p$,
for every prime $p \neq p_0$, then $\mathcal{C}$ does
\textbf{not} split over~$\real$. But we have no need for this
more general (and somewhat less concise) version.
 \end{rem}

\section{The method of Pink and Prasad}
\label{PinkPrasadMethSect}

\begin{proof}[{\bf Proof of Prop.~\ref{Quasi->good}}]
  Suppose $(\pi, V_\rational)$ is an irreducible
representation of~$G$ over~$\rational$. Let 
 $$ \mathcal{C} = \End_{G}(V_\rational) $$
 be the centralizer of $\pi(G_\rational)$ in
$\End_\rational(V_\rational)$. Then Schur's Lemma tells us that
$\mathcal{C}$~is a division algebra over~$\rational$.

Because $G$ splits over the quadratic extension~$F$,
we know that $V_\complex$ is either irreducible or the sum of two
irreducibles \fullsee{Q(i)split}{2irreds}. 

\setcounter{case}{0}

\begin{case}
 Assume $V_\complex$ is irreducible.
 \end{case}
 Then $V_\real$ is obviously irreducible.

\begin{case}
 Assume $V_\complex$ is the direct sum of two irreducibles that
are not isomorphic. 
 \end{case}
 From the assumption of this case, and the fact that $G$~splits
over~$F$, we know that $V_F$ is the direct sum of two
irreducibles that are not isomorphic. Therefore, $\End_G(V_F)
\iso F \oplus F$, so $\mathcal{C} \otimesQ F \iso F \oplus F$. 

Write $F = \rational[\sqrt{-r}]$, for some $r \in \rational^+$. 
 \begin{itemize}
 \item Because $F \oplus F$ is commutative, we know that
$\mathcal{C}$~is commutative, so $\mathcal{C}$~is a field.
 \item Because $\dim_F (F \oplus F) = 2$, we know that
$\dim_\rational \mathcal{C} = 2$. 
 \item Because $\mathcal{C} \otimesQ F = \mathcal{C}[\sqrt{-r}]$
is not a field, we know that $\mathcal{C}$~contains a root of
$x^2 + r$. 
 \end{itemize}
 We conclude that $\mathcal{C} \iso F$.

Therefore 
 $$ \End_G(V_\real) \iso \mathcal{C} \otimesQ \real \iso F
\otimesQ \real \iso \complex $$
 is a field. So $V_\real$ is irreducible.

\begin{case}
 Assume $V_\complex$ is the direct sum of two irreducibles that
are isomorphic. 
 \end{case}
  In this case, we know that $\End_G(V_\complex) \iso \Mat_{2
\times 2}(\complex)$ is $4$-dimensional over~$\complex$. Since
$\End_G(V_\complex) \iso \mathcal{C} \otimesQ \complex$, we
conclude that $\mathcal{C}$ is $4$-dimensional over~$\rational$.
Thus, $\mathcal{C}$ is a quaternion algebra over~$\rational$.

For every odd prime~$p$, Prop.~\ref{QuasiSplit->Field} (and the
fact that $\mathcal{C} \otimesQ \rational_p$ is not commutative)
implies that $V_{\rational_p}$ is reducible, so $\mathcal{C}
\otimesQ \rational_p = \End_G(V_{\rational_p})$ is not a
division algebra. In other words, $\mathcal{C}$ splits
over~$\rational_p$. 

Now, Lem.~\ref{SplitLocal->Split} asserts that $\mathcal{C}$
does not split over~$\real$. This means that $\End_G(V_\real)
\iso \mathcal{C} \otimesQ \real$ is a division algebra. We
conclude that $V_\real$ is irreducible, as desired.
 \end{proof}

\begin{rem}
 From the proof (and Rem.~\ref{HasseAnyPrime}), it is clear that
the exceptional prime~$2$ can be replaced with any other prime in
Condition~\fullref{Quasi->good}{quasi}: it suffices to assume
that there is a prime~$p_0$, such that $G$ is quasi-split over
the $p$-adic field $\rational_p$, for every prime $p \neq p_0$.
The case $p_0 = 2$ is all we need for our proof of
Cor.~\ref{EGoodQform}.
 \end{rem}

\section{Reducing to the compact case}
\label{RaghuMethodSect}

\begin{defn} \label{LongWeylDefn}
 Suppose $G$ is a connected, reductive algebraic
$\rational$-group. Let 
 \begin{itemize} 
 \item $S$~be a maximal $\rational$-split torus of~$G$;
 \item $C = C_G(S)$ be the centralizer of~$S$;
 \item $M'$ be the (unique) maximal connected, semisimple
subgroup of the reductive group~$C$;
 \item $T$~be a maximal $\rational$-torus of~$M'$;
 \item $\Phi^+$ be the positive roots of $(\Lie M'_\complex,
\Lie T_\complex)$ (with respect to some ordering);
 and 
 \item $\Phi^-$ be the set of negative roots.
 \end{itemize}

We call $M'$ the \emph{semisimple anisotropic kernel} of~$G$.

 We say that \emph{the longest element of the Weyl group of
the anisotropic kernel of~$G$ is realized over~$\rational$} if
there is some $w \in N_{M'}(T)_\rational$, such that 
 $$ w(\Phi^+) = \Phi^- .$$
 Here, as usual, the normalizer $N_{M'}(T)$ acts on $\Lie T^*$
by 
 $$w(\lambda)(t) = \lambda(w^{-1} t w) .$$
 It is important to
notice (from the subscripts in $N_{M'}(T)_\rational$) that
$w$~is required to be in the semisimple group~$M'$, and that
$w$~is required to be a $\rational$-element.
 \end{defn}

\begin{thm} \label{Weyl->good}
 Suppose $G$ is a connected, reductive algebraic
$\rational$-group. If
 \begin{enumerate} \renewcommand{\theenumi}{\alph{enumi}}
 \item $G$ is split over some imaginary quadratic
extension~$F$ of~$\rational$,
 \item $\Qrank G = \Rrank G$,
 and
 \item the longest element of the Weyl group of
the anisotropic kernel of~$G$ is realized over~$\rational$,
 \end{enumerate}  
 then each irreducible $\rational$-representation of~$G$
remains irreducible over~$\real$.
 \end{thm}

M.~S.~Raghunathan \cite[\S3]{Raghunathan-Qform} proved
Thm.~\ref{Weyl->good} in the special case where 
 \begin{itemize}
 \item $G$ is semisimple,
 and
 \item $\Rrank G = 0$ (in other words, $G$~is compact).
 \end{itemize}
 The following result shows that the general case follows from
this.

\begin{prop} \label{AnisKer->G}
 Suppose $G$ is a connected, reductive algebraic group
over~$\rational$, and let $M'$ be the semisimple anisotropic
kernel of~$G$. If
 \begin{itemize}
 \item $G$ is split over some imaginary quadratic extension~$F$
of~$\rational$;
 \item $\Qrank G = \Rrank G$;
 and
 \item every irreducible $\rational$-representation of~$M'$
remains irreducible over~$\real$,
 \end{itemize}
 then every irreducible $\rational$-representation of~$G$
remains irreducible over~$\real$.
 \end{prop}

\begin{rem}
 There is no need to assume that the quadratic extension~$F$ is
imaginary in \pref{Weyl->good} or~\pref{AnisKer->G}: if $F$~is
real, then the hypotheses imply that $G$~is $\rational$-split,
so Lem.~\ref{split->reps} applies.
 \end{rem}

Before proving the proposition, let us state a simple
lemma, which reduces the construction of $\rational$-forms of
representations of~$G$ to the same problem for certain
representations of a minimal parabolic~$P$. It is similar to the
usual construction of highest weight modules.

\begin{lem} \label{HighUt->Qform}
 Let 
 \begin{itemize}
 \item $\Lie G$ be a semisimple Lie algebra over~$\rational$;
 \item $\Lie T$ be a maximal $\rational$-split torus of~$\Lie
G$;
 \item $\Phi_\rational$ be the system of $\rational$-roots of
$(\Lie G, \Lie T)$;
 \item $\Lie P$~be a minimal parabolic $\rational$-subalgebra
of~$\Lie G$ that contains~$\Lie T$;
 \item $V$ be an irreducible, real $\Lie G$-module;
 \item $\lambda \colon \Lie T \to \rational$ be the highest
weight of~$V$, with respect to the ordering
of~$\Phi_\rational$ determined by the parabolic~$\Lie P$;
 and
 \item $U$ be a $\Lie P$-invariant $\rational$-form of the
weight space~$V^\lambda$.
 \end{itemize}
 Then the representation~$V$ has a $\rational$-form.
 \end{lem}

\begin{proof}
 Let $\Delta_\rational$ be the base of~$\Phi_\rational$
determined by~$\Lie P$. Let
 $\widetilde U$ be the $\rational$-span of 
 $$ \{\, y_1 y_2 \cdots y_k w \mid w \in U, \alpha_j \in
\Delta_\rational, y_j \in \Lie G_{-\alpha_j} \,\}
.$$

\setcounter{step}{0}

\begin{step} \label{UisGinv}
  $\widetilde U$ is $\Lie G$-invariant.
 \end{step}
 From the definition of~$\widetilde U$, it is obvious that
$[\Lie G_{-\alpha}, \widetilde U] \subset \widetilde U$ for
all $\alpha \in \Delta_\rational$. Also, because both the
centralizer ${\Lie C}_{\Lie G}(\Lie T)$ and the root space~$\Lie
G_\alpha$ are contained in~$\Lie P$, it is not difficult to see
(by induction on~$k$) that
 $[\Lie C_{\Lie G}(\Lie T), \widetilde U] \subset \widetilde
U$ and
 $[\Lie G_\alpha, \widetilde U] \subset \widetilde U$, for all
$\alpha \in \Delta_\rational$. Since 
 $$ \mbox{$\Lie C_{\Lie G}(\Lie T) \cup \bigcup_{\alpha \in
\Delta_\rational} (\Lie G_\alpha \cup \Lie G_{-\alpha})$
generates $\Lie G$,} $$
 we conclude that $\widetilde U$ is $\Lie G$-invariant.

\begin{step} \label{UspansV}
 $\widetilde U$ spans~$V$ over~$\real$.
 \end{step}
 The $\real$-span of~$\widetilde U$ is a submodule of~$V$. so
the desired conclusion follows from the fact that $V$ is
irreducible.

\begin{step} \label{Uindep/Q}
 If $a_1,\ldots,a_r$ are real numbers that are linearly
independent over~$\rational$, and $w_1,\ldots,w_r$ are
nonzero elements of~$\widetilde U$, then $\sum_{j=1}^r a_j
w_j \neq 0$.
 \end{step}
 Suppose $\sum_{j=1}^r a_j w_j = 0$. (This will lead to a
contradiction.)

 Since 
 $y_1 y_2 \cdots y_k w \in V^{\lambda - \alpha_1 - \cdots -
\alpha_k}$, we see that 
 \begin{equation} \label{UcapHigh=U}
 \widetilde U \cap V^\lambda = U 
 \end{equation}
 and
 \begin{equation} \label{UisGraded}
 \widetilde U = \bigoplus_{\mu \in \Lie T^*} (\widetilde U
\cap V^\mu) .
 \end{equation}
 Because of \pref{UisGraded}, we may assume there is some
weight~$\mu$, such that $w_j \in V^\mu$ for all~$j$. (Project
to some~$V^\mu$, and delete the $w_j$'s whose projection
is~$0$.) 

Because $V$~is irreducible, and $\lambda$~is the highest
weight, there exist $x_1,\ldots,x_k \in \cup_{\alpha \in
\Delta_\rational} \Lie G_\alpha$, such that $x_1 \cdots x_k
w_1$ is a nonzero element of~$V^\lambda$. From
\pref{UisGraded}, we see that $x_1 \cdots x_k w_j \in
V^\lambda$ for every~$j$. Hence, \pref{UcapHigh=U} implies
that $x_1 \cdots x_k w_j \in U$ for all~$j$. Since
 $$ \sum_{j=1}^r a_j (x_1 \cdots x_k w_j)
 = x_1 \cdots x_k \sum_{j=1}^r a_j w_j
 =  x_1 \cdots x_k \cdot 0
 = 0 , $$
 and $U$ is a $\rational$-form of~$V^\lambda$, this implies
that $x_1 \cdots x_k w_j = 0$ for every~$j$. This contradicts
the choice of $x_1,\ldots,x_k$.

\begin{step}
 Completion of proof.
 \end{step}
 From Steps~\ref{UspansV} and~\ref{Uindep/Q}, we see that the
natural scalar-multiplication map $\real \otimesQ
\widetilde U \to V$ is a bijection. So $\widetilde U$ is a
$\rational$-form of the vector space~$V$. By combining this
with Step~\ref{UisGinv}, we conclude that $\widetilde U$ is a
$\rational$-form of the representation.
 \end{proof}

\begin{proof}[{\bf Proof of Prop.~\ref{AnisKer->G}}]
 We consider two cases.

\setcounter{case}{0}

\begin{case} \label{AnisKer->G-cpct}
 Assume that the semisimple part of~$G_\real$ is compact.
 \end{case}
 We may write $G = M'  A T$, where 
 \begin{itemize}
 \item $M'$~is connected and semisimple (so, by assumption,
$M'_\real$~is compact);
 \item $A$~is a $\rational$-split torus;
 and
 \item $T$~is a torus that is anisotropic over~$\rational$.
 \end{itemize}
 Let $V$ be an irreducible $\rational$-representation of~$G$.
Since $A$, being a $\rational$-split, central torus, acts by
scalars, we know that $V$ is an irreducible
$\rational$-representation of $M' T$.

\begin{subcase}
 Assume $T$ acts trivially on~$V$. 
 \end{subcase}
 Then $V$ is an irreducible $\rational$-representation of~$M'$,
so, by assumption, it remains irreducible over~$\real$.

\begin{subcase}
 Assume $T$ acts nontrivially on~$V$. 
 \end{subcase}
 Because $T_F$ is an $F$-split, central torus, its Lie algebra
defines an action of~$F$ on~$V$ that centralizes~$M'$. Thus,
we may think of~$V$ as an irreducible $F$-representation of~$M'$.
Let us say $V_\rational = W_F$; then $V_\real = W_\complex$.
Because $M'$ is $F$-split, we know that $W_\complex$ is
irreducible \fullsee{Q(i)split}{irred}.

Because $\Qrank T = 0$ and $\Qrank G = \Rrank G$, we see that
$\Rrank T = 0$, which means $T_\real$~is compact, so the Lie
algebra of~$T$ acts by purely imaginary scalars on~$W_\complex$.
Thus, any $M' T$-invariant $\real$-submodule of $V_\real$ is an
$M'$-invariant $\complex$-submodule of~$W_\complex$. Hence, the
conclusion of the preceding paragraph implies that $V_\real$ is
irreducible.

\begin{case}
 The general case.
 \end{case}
 Given a representation~$V$ of~$G$ over~$\real$,
we wish to show that $V$ has a $\rational$-form (see
Lem.~\ref{Qform<>RIrred}). Because representations of~$G$ are
completely reducible, we may assume that $V$~is irreducible.

Let $P$ be a minimal parabolic $\rational$-subgroup of~$G$,
let $T$~be a maximal $\rational$-split torus of~$P$, and
$\lambda$ be the highest weight of~$V$ (with respect to~$T$
and~$P$).

Now $V^\lambda$ is $C_G(T)$-invariant, so, from
Case~\ref{AnisKer->G-cpct}, we know that the vector
space~$V^\lambda$ has a $\rational$-form~$U$ that is
$C_G(T)_\rational$-invariant. Then, since the unipotent
radical of~$P$ annihilates~$V^\lambda$, we know that $U$ is
$P_\rational$-invariant. So Lem.~\ref{HighUt->Qform} implies
that $V$~ has a $\rational$-form.
 \end{proof}

\section{Construction of a good $\rational$-form}
\label{ConstructQformSect}

In this section, we provide an explicit construction of a
$\rational$-form of~$G$ that satisfies the hypotheses of
Prop.~\ref{Quasi->good} and Thm.~\ref{Weyl->good}. 

\begin{prop} \label{ConstructQform}
 If $G$ is a connected, simply connected, semisimple algebraic
$\real$-group, then $G$ has a $\rational$-form, such that
 \begin{enumerate}
 \item \label{ConstructQform-Qisplit}
 $G$ is split over $\Qi$,
 \item \label{ConstructQform-quasi}
 $G$ is quasi-split over the $p$-adic field $\rational_p$,
for every odd prime~$p$,
  \item \label{ConstructQform-Qrank}
 $\Qrank G = \Rrank G$,
 and
 \item \label{ConstructQform-Weyl}
 every element of the Weyl group of the anisotropic kernel
of~$G$ is defined over~$\rational$.
 \end{enumerate}
 \end{prop}

The argument is a straightforward adaptation of A.~Borel's
\cite{Borel-CK} classical proof of the existence of an
anisotropic $\rational$-form.

Actually, like Borel, we do not directly construct
$G_\rational$ itself, but only the Lie algebra~$\Lie
G_\rational$, so, to avoid problems in passing from the Lie
algebra to the group, we need to make some assumption on the
fundamental group of~$G$. (It needs to be a
$\rational$-subgroup of the universal cover~$\widetilde G$.)
Therefore, the statement of Prop.~\ref{ConstructQform} requires
$G$ to be simply connected. Alternatively, one could
require~$G$ to be adjoint, instead of simply connected, but
the situation is not obvious for some intermediate groups
that are neither adjoint nor simply connected.

\begin{rem}
 As a complement to our explicit construction, it might be
possible to use theorems of Galois cohomology to give a more
elegant proof of Prop.~\ref{GRhasQform}. In this vein, G.~Prasad
(see \cite[Prop.~6.4]{OhAdelic}) gave a very short proof of the
existence of a $\rational$-form satisfying
\fullref{ConstructQform}{Qisplit}
and~\fullref{ConstructQform}{Qrank}; perhaps a clever argument
can yield \fullref{ConstructQform}{quasi}
and/or~\fullref{ConstructQform}{Weyl}, as well, but they do not
seem to be obvious.
 \end{rem}

Let us set up the usual notation.

\begin{notation} \label{SSnotation}
 \ 
 \begin{itemize}
 \item $\Lie G$ is a real semisimple Lie algebra;
 \item $\kappa(\cdot,\cdot)$ is the Killing form on~$\Lie G$;
 \item $\Lie H$ is a maximal torus (i.e., a Cartan
subalgebra) of~$\Lie G$;
 \item $\Phi$ is the set of roots of $(\Lie G_\complex, \Lie
H_\complex)$;
 \item $h_\alpha$ is the unique element of~$\Lie H_\complex$,
such that $\alpha(t) = \kappa(t, h_\alpha)$ for all $t \in
\Lie H_\complex$ (for each $\alpha \in \Phi$);
 \item $h_\alpha^* = 2 h_\alpha/\kappa(h_\alpha,h_\alpha)$
(for each $\alpha \in \Phi$);
 \item $(\Lie G_\complex)_\alpha$ is the root space
corresponding to $\alpha \in \Phi$;
 \item $\theta$ is a Cartan involution of~$\Lie G$, such that
$\theta(\Lie H) = \Lie H$ (we also use $\theta$~to denote the
extension to a $\complex$-linear automorphism of~$\Lie
G_\complex$);
 \item $\Lie G = \Lie K + \Lie P$ is the Cartan decomposition
of~$\Lie G$ corresponding to~$\theta$ (i.e., $\Lie K$
and~$\Lie P$ are the $+1$ and~$-1$ eigenspaces of~$\theta$,
respectively).
 \end{itemize}
 Because $\theta(\Lie H) = \Lie H$, we see that
$\theta$~induces a permutation of~$\Phi$: we have
 $\theta \bigl( (\Lie G_\complex)_\alpha \bigr)
 = (\Lie G_\complex)_{\theta(\alpha)}$.
 \end{notation}

The following lemma is a slight modification of a result of
Borel \cite[\S3.2 and Lem.~3.5]{Borel-CK} that extends work
of Chevalley and Weyl. (See \cite[p.~116 and footnote on
p.~117]{Borel-CK} for some historical remarks.) We follow
Borel's proof almost verbatim. However, Borel assumed that
the Cartan subalgebra~$\Lie H$ contains a maximal
$\real$-anisotropic torus of~$\Lie G$, and, using this
assumption, he obtained a stronger version
of~\pref{Borelbasis-theta}: $\theta(x_\alpha) = \pm
x_{\theta(\alpha)}$.

\begin{lem}[Borel, Chevalley, Gantmacher, Weyl]
 \label{Borelbasis}
 Assume the notation of \pref{SSnotation}.

 There is a function
 $\Phi \to \Lie G \colon \alpha \mapsto x_\alpha$, such that,
for $\alpha,\beta \in \Phi$, we have
 \begin{enumerate}
 \item \label{Borelbasis-xisroot}
 $x_\alpha \in \Lie G_\alpha$;
 \item \label{Borelbasis-brackets}
 $ \displaystyle
 [x_\alpha,x_\beta] = 
 \begin{cases}
 N_{\alpha,\beta} x_{\alpha+\beta}
 & \text{if $\alpha + \beta \in \Phi$} , \\
 \hfil -h_\alpha^*
 & \text{if $\alpha + \beta = 0$} , \\
 \hfil 0
 & \text{if $0 \neq \alpha + \beta \notin \Phi$} ,
 \end{cases}
 $
 \newline
 where
 \begin{equation}
 N_{\alpha,\beta} = N_{-\alpha,-\beta} =
\pm(p_{\alpha,\beta}+1)
 \end{equation}
 and $p_{\alpha,\beta} \ge 0$ is the greatest integer such
that $\alpha - p_{\alpha,\beta} \beta \in \Phi$;
 \smallskip
 \item \label{Borelbasis-theta}
 $\theta(x_\alpha) \in \{ \pm x_{\theta(\alpha)}, \pm i
x_{\theta(\alpha)} \}$;
 and
 \smallskip
 \item \label{Borelbasis-K+iP}
 $ \displaystyle 
 \Lie K + i\Lie P
 = \sum_{\alpha \in \Phi} i \real h_\alpha
 + \sum_{\alpha \in \Phi}
 \bigl\{\, z x_\alpha + \conj{z} x_{-\alpha} \mid z \in
\complex \, \bigr\}
 $.
 \end{enumerate}
 \end{lem}

\begin{proof}[Proof {\upshape(\cite[\S3.2--\S3.5]{Borel-CK}
or \cite[Chap.~14]{RaghunathanBook})}]
 The famous Chevalley basis \cite{ChevalleyTohoku} satisfies
\pref{Borelbasis-xisroot} and~\pref{Borelbasis-brackets}. 

\setcounter{step}{0}

\begin{step}
 We may assume \pref{Borelbasis-K+iP} holds.
 \end{step}
 Recall that all of the maximal compact subgroups of any
connected Lie group are conjugate to each other, and that all
of the maximal toruses of any connected, compact Lie group
are conjugate to each other. Thus, since the LHS and RHS of
\pref{Borelbasis-K+iP} are maximal compact subalgebras
of~$\Lie G_\complex$ that contain the maximal torus
$\sum_{\alpha \in \Phi} i \real h_\alpha$, they are
conjugate, via an automorphism of $\Lie G_\complex$ that
normalizes $\sum_{\alpha \in \Phi} i \real h_\alpha$. Hence,
by replacing $\{x_\alpha\}_{\alpha \in \Phi}$ with a
conjugate, we may assume \pref{Borelbasis-K+iP} holds.

\begin{step}
 For each $\alpha \in \Phi$, define $c_\alpha \in \complex$ by
$\theta(x_\alpha) = c_\alpha x_{\theta(\alpha)}$; then
 \begin{equation} \label{c(-a)c(a)=1}
  \mbox{$\displaystyle
 c_{-\alpha} = \frac{1}{c_\alpha} = c_{\theta(\alpha)}$
 \  for all $\alpha \in \Phi$}
 \end{equation}
 and
 \begin{equation} \label{c(a+b)}
   \mbox{$\displaystyle
 c_\alpha c_\beta = \pm c_{\alpha + \beta}$
 \  for all $\alpha,\beta \in \Phi$, such that $\alpha+\beta
\in \Phi$} .
 \end{equation}
 \end{step}
 Note that
 \begin{align*}
 - c_\alpha c_{-\alpha} h_{\theta(\alpha)}^*
 &= [c_\alpha x_{\theta(\alpha)}, c_{-\alpha} 
x_{-\theta(\alpha)}] \\
 &= [\theta(x_\alpha), \theta(x_{-\alpha})]
 = \theta [x_\alpha,x_{-\alpha}]
 = - \theta(h_\alpha^*)
  . 
 \end{align*}
 Because $\theta(h_\alpha^*) = h_{\theta(\alpha)}^*$ (since
$\theta$ is an automorphism that fixes~$\Lie H$), this implies
that $c_\alpha c_{-\alpha} = 1$, which establishes part of
\pref{c(-a)c(a)=1}. For the other part, we use the
fact that $\theta^2 = \Id$ to calculate
 $$ x_\alpha = \theta^2(x_\alpha) = c_\alpha \,
\theta(x_{\theta(\alpha)})
 = c_\alpha \, c_{\theta(\alpha)} \, x_\alpha .$$

To establish \pref{c(a+b)}, note that
$p_{\theta(\alpha),\theta(\beta)} = p_{\alpha,\beta}$
(because $\theta$~is an automorphism), so
$N_{\theta(\alpha),\theta(\beta)} = \pm N_{\alpha,\beta}$.
Now use the fact that
 $$ [\theta(x_\alpha), \theta(x_\beta)] 
 = \theta[x_\alpha,x_\beta] .$$

\begin{step}
 We may assume \pref{Borelbasis-theta} holds.
 \end{step}
 Let $\Delta$ be a basis of~$\Phi$ (with respect to some
order). Then $\Delta$ is a basis of the dual space~$\Lie
H^*$, so there is some $h \in \Lie h_{\complex}$, such
that $e^{\alpha(h)} = c_\alpha$, for every $\alpha \in
\Delta$. Then, from~\pref{c(a+b)}, we see, by induction on
the length of~$\alpha$, that $e^{\alpha(h)} = \pm
c_\alpha$, for every $\alpha \in \Phi^+$. Because
$c_{-\alpha} = 1/c_\alpha$, then we have
  $e^{\alpha(h)} = \pm
c_\alpha$, for every $\alpha \in \Phi$.

 For each $\alpha \in \Phi$, let $x'_\alpha =
e^{-\alpha(h)/2} x_\alpha$. Then it is easy to see that
\pref{Borelbasis-xisroot} and~\pref{Borelbasis-brackets}
hold with $x'_\alpha$ in the place of~$x_\alpha$. 

 Because
 $\theta(\Lie K + i \Lie P) = \Lie K + i \Lie P$, we see
from~\pref{Borelbasis-K+iP} that $c_{-\alpha} =
\conj{c_\alpha}$.
 Then
  $1/c_\alpha = c_{-\alpha} = \conj{c_\alpha}$, so
$|c_\alpha| = 1$. Therefore, $\alpha(h)$ is pure imaginary
for every $\alpha \in \Delta$. By linearity, then $\alpha(h)$
is pure imaginary for every $\alpha \in \Phi$, so
$e^{-\alpha(h)/2} = \conj{e^{\alpha(h)/2}}$, for every
$\alpha \in \Phi$. Thus, \pref{Borelbasis-K+iP} holds with
$x'_\alpha$ in the place of~$x_\alpha$. 

For any $\alpha \in \Phi$, we have
 \begin{align*}
 \theta(x'_\alpha)
 &= e^{-\alpha(h)/2} \, \theta(x_\alpha) \\
 &= (\pm c_\alpha)^{-1/2} \bigl( c_\alpha x_{\theta(\alpha)}
\bigr) \\
 &= (\pm c_\alpha)^{1/2}
 \bigl( e^{\theta(\alpha)(h)/2} x'_{\theta(\alpha)} \bigr) \\
 &= (\pm c_\alpha)^{1/2}
 ( \pm c_{\theta(\alpha)})^{1/2} x'_{\theta(\alpha)} \\
 & =  (\pm c_\alpha)^{1/2}
 (\pm c_\alpha)^{-1/2} x'_{\theta(\alpha)} \\
 &=  (\pm 1)^{1/2}  x'_{\theta(\alpha)} .
 \end{align*}
 Thus, \pref{Borelbasis-theta} holds with $x'_\alpha$ in the
place of~$x_\alpha$.
 \end{proof}

Prop.~\ref{ConstructQform} is obtained quite easily from this
lemma. Most of the argument we give is based on \cite{Borel-CK}
or \cite[Chap.~14]{RaghunathanBook}. However, Steps
\ref{LittleSU2} and~\ref{anisWeyl} are from
\cite[\S2]{Raghunathan-Qform}, and Steps~\ref{LittleSplit}
and~\ref{Qp-quasisplit} are based on suggestions of G.~Prasad
(personal communication).

\begin{proof}[{\bf Proof of Prop.~\ref{ConstructQform}}]
 We begin by establishing notation.
 \begin{itemize}
 \item Let $\Lie G$ be the Lie algebra of~$G$.
 \item Choose a maximal torus~$\Lie H$ of~$\Lie G$, such that
$\Rrank \Lie H = \Rrank \Lie G$.
 \item Assume the notation of~\pref{SSnotation}.
 \item Let $\Lie T = \Lie H \cap \Lie P$ (so $\Lie T$ is a
maximal $\real$-split torus of~$\Lie G$).
 \item Let $\{x_\alpha\}_{\alpha \in \Phi}$ be as in
Lem.~\ref{Borelbasis}.
 \item Let $\Lie G_\Qi$ be the $\Qi$-span of $\{h_\alpha^*,
x_\alpha\}_{\alpha \in \Phi}$ in $\Lie G_\complex = \Lie G
\otimesR \complex$.
 \item Let $\Lie H_\Qi = \Lie H_\complex \cap \Lie G_\Qi$ be
the $\Qi$-span of $\{h_\alpha^*\}_{\alpha \in \Phi}$ in $\Lie
G_\complex$.
 \item Let $\Lie G_\rational = \Lie G_\Qi \cap \Lie G$.
 \end{itemize}
 We will show that $\Lie G_\rational$ is a $\rational$-form
of~$\Lie G$, and that the corresponding $\rational$-form
$G_\rational$ of~$G$ satisfies all the hypotheses of
Prop.~\ref{ConstructQform}.

\setcounter{step}{0}

\begin{step} \label{GQisplit}
 $\Lie G_\Qi$ is a split $\Qi$-form of~$\Lie G_\complex$,
with $\Lie H_\Qi$~being a $\Qi$-split maximal torus.
 \end{step}
 It is clear that $\Lie G_\Qi$ is a $\Qi$-form of~$\Lie
G_\complex$ (because it is the $\Qi$-span of a basis, and is
closed under brackets). It is split because it contains a
Chevalley basis of~$\Lie G_\complex$, and $\Lie H_\Qi$~is the
maximal split torus corresponding to this basis.

\begin{step} \label{K&P/Q}
 Each of $\Lie K$ and $\Lie P$ is the  $\real$-span of its intersection
with~$\Lie G_\rational$.
 \end{step}
 Because each of these subspaces is contained in~$\Lie G$, it
suffices to prove the conclusion with $\Lie G_\rational$
replaced by $\Lie G_\Qi$.

 Let $U_\rational = (\Lie K + i \Lie P) \cap \Lie G_\Qi$.
Then $U_\rational$~is a $\rational$-subspace of~$\Lie K + i
\Lie P$. Indeed, we see, from \fullref{Borelbasis}{K+iP},
that $U_\rational$~is a $\rational$-form of the real vector
space $\Lie K + i \Lie P$.

Now $\Lie K + i \Lie P$ and~$\Lie G_\Qi$ are
$\theta$-invariant (for the latter, see
\fullref{Borelbasis}{theta}). So $U_\rational$~is
$\theta$-invariant. This means that, with respect to the
$\rational$-form~$U_\rational$, the linear transformation
$\theta|_{\Lie K + i \Lie P}$ is defined over~$\rational$.
Since the eigenvalues ($\pm 1$) are rational, we conclude
that the eigenspaces are spanned (over~$\real$) by the
rational vectors, that is by elements of~$U_\rational$.
Concretely, this means that the $\real$-span of $\Lie K \cap
U_\rational$ is~$\Lie K$, and the $\real$-span of $i\Lie P
\cap U_\rational$ is~$i \Lie P$. The first is exactly what we
want to know about~$\Lie K$. Multiplying by~$i$ transforms
the second into exactly what we want to know about~$\Lie P$.

\begin{step} \label{GQ=Qform}
 $\Lie G_\rational$ is a $\rational$-form of~$\Lie G$.
 \end{step}
 Because $\Lie G$ and $\Lie G_\Qi$ are closed under brackets,
it is clear that $\Lie G_\rational$ is a subalgebra of~$\Lie
G$. We just need to show that its $\real$-span is all
of~$\Lie G$.

From Step~\ref{K&P/Q}, we know that the $\real$-span of $\Lie
G_\rational$ contains both~$\Lie K$ and~$\Lie P$. Therefore,
it contains $\Lie K + \Lie P = \Lie G$.

\begin{step}
 $\Lie G_\rational$ splits over~$\Qi$.
 \end{step}
 We already pointed out in Step~\ref{GQisplit} that $\Lie
G_\Qi$ is split.

\begin{step}
 $\Qrank \Lie G_\rational = \Rrank \Lie G$.
 \end{step}
 Because $\Lie T = \Lie H \cap \Lie P$ is a maximal
$\real$-split torus of~$\Lie G$, it suffices to show that
$\Lie T$ is (defined over~$\rational$ and) $\rational$-split.

\begin{substep}
 $\Lie T$ is defined over~$\rational$.
 \end{substep}
 From Step~\ref{GQ=Qform}, we see that $\Lie G_\Qi = \Lie
G_\rational + i \Lie G_\rational$, so, for any (real)
subspace~$X$ of~$\Lie G$, we have $X_\complex \cap \Lie G_\Qi
= (X \cap \Lie G_\rational)_\complex$. Thus, if $X_\complex$
is the $\real$-span of its intersection with $\Lie G_\Qi$,
then $X$~is the $\real$-span of its intersection with $\Lie
G_\rational$, i.e., $X$~is defined over~$\rational$.

It is clear, from the definition of~$\Lie G_\Qi$, that $\Lie
H_\complex$ is the $\real$-span of its intersection with
$\Lie G_\Qi$. Hence, from the preceding paragraph, we
conclude that $\Lie H$~is defined over~$\rational$.
 From Step~\ref{K&P/Q}, we know that $\Lie P$ is also defined
over~$\rational$. Hence, the intersection $\Lie T = \Lie H
\cap \Lie P$ is defined over~$\rational$.

\begin{substep}
 $\Lie T$ is is $\rational$-split.
 \end{substep}
 Let $T$ be the $\rational$-torus of~$G$ corresponding
to~$\Lie T$. We know that $T$ splits over~$\real$, so
$\chi(T_\real) \subset \real$, for every character~$\chi$
of~$T$. Because $\Lie T \subset \Lie H$, and $\Lie H$~splits
over $\Qi$ (see Step~\ref{GQisplit}), we know that
$\chi(T_\rational) \subset \Qi$, for every character~$\chi$
of~$T$. So $\chi(T_\rational) \subset \real \cap \Qi =
\rational$, for every character~$\chi$ of~$T$; hence, $T$~is
$\rational$-split.

\begin{step} \label{LittleSU2}
 For each $\alpha \in \Phi$, let
 $$
  \Lie G_{\pm \alpha} = \langle (\Lie G_\complex)_\alpha,
(\Lie G_\complex)_{-\alpha} \rangle \cap \Lie G {.}
 $$
 If
$\alpha(\Lie T) = 0$, then $\bigl( \Lie G_{\pm \alpha}
\bigr)_\rational \iso \su(2)_\rational$
  {\upshape(}for the usual $\rational$-form on
$\SU(2)${\upshape)}.
 \end{step}
 Because $\Lie T$ is a maximal $\real$-split torus, the
assumption on~$\alpha$ implies that $(\Lie G_\complex)_\alpha
\subset \Lie K_\complex$ (and the same for~$-\alpha$). So,
using \fullref{Borelbasis}{K+iP}, we see that
 \begin{align*}
 \langle (\Lie G_\complex)_\alpha, (\Lie
G_\complex)_{-\alpha} \rangle \cap \Lie G
 &= \langle (\Lie G_\complex)_\alpha, (\Lie
G_\complex)_{-\alpha} \rangle \cap \Lie K \\
 &= i\real h_\alpha^* + \{\, z x_\alpha + \conj{z} x_{-\alpha}
\mid z \in \complex \,\} .
 \end{align*}
 Therefore
 \begin{align*}
 \bigl( \langle (\Lie G_\complex)_\alpha, (\Lie
G_\complex)_{-\alpha} \rangle \cap \Lie G \bigr)_\rational
 &= i\rational h_\alpha^* + \{\, z x_\alpha + \conj{z}
x_{-\alpha} \mid z \in \Qi \,\} \\
 &\iso \su(2)_\rational .
 \end{align*}

\begin{step} \label{anisWeyl}
 Every element of the Weyl group of the anisotropic kernel
of~$G$ is realized over~$\rational$.
 \end{step}
 The Weyl element of $\SU(2)$ is realized by the rational
matrix
 $$ \begin{bmatrix}
 0 & 1 \\
 -1 & 0 
 \end{bmatrix}
 .$$
 So Step~\ref{LittleSU2} implies that all of the root
reflections of the anisotropic kernel can be realized
over~$\rational$. These reflections generate the entire Weyl
group.

\begin{step} \label{LittleSplit}
 For any odd prime~$p$, and any $\alpha \in \Phi$, such that
$\alpha(\Lie T) = 0$, the Lie algebra $\Lie G_{\pm \alpha}$ is
$\rational_p$-split.
 \end{step}
 The quadratic form $x_1^2 + x_2^2 + x_3^2$ is isotropic
over~$\rational_p$ (see, for example, \cite[Cor.~1.6.2,
p.~50]{BorevichShafarevich}), so $\so(3)$ is
$\rational_p$-split. Since $\so(3)_{\rational_p} \iso
\su(2)_{\rational_p}$ and $\su(2)_{\rational_p} \iso (\Lie G_{\pm
\alpha})_{\rational_p}$ (see Step~\ref{LittleSU2}), we conclude
that $\Lie G_{\pm \alpha}$ is $\rational_p$-split.

\begin{step} \label{Qp-quasisplit}
 $\Lie G_{\rational_p}$ is quasi-split, for every odd prime~$p$.
 \end{step}
 Let $\Psi$ be a maximal set of pairwise orthogonal roots in 
 $\{\, \alpha \in \Phi \mid \alpha(\Lie T) = 0 \,\}$.
 For each $\alpha \in \Psi$, we know, from
Step~\ref{LittleSplit}, that $(\Lie G_{\pm
\alpha})_{\rational_p}$ contains a nontrivial
$\rational_p$-split torus $\Lie S_\alpha$; let 
 $$\Lie S_{\rational_p} = \Lie T_{\rational_p} + \sum_{\alpha
\in \Psi} \Lie S_\alpha .$$
 From the maximality of~$\Psi$, we know that the centralizer of
the torus 
 $$\Lie S' = \Lie T_{\rational_p} + \sum_{\alpha \in \Psi} i
\rational_p h^*_\alpha$$
 in $\Lie G_{\rational_p}$ is $\Lie H_{\rational_p}$, a (maximal)
torus of~$\Lie G_{\rational_p}$. Now both $\Lie S$ and $\Lie S'$
are maximal tori of the Lie algebra
 $$\Lie T_{\rational_p} +  \sum_{\alpha \in \Psi} \Lie
G_{\pm\alpha} ,$$
 so they are conjugate over the algebraic closure
of~$\rational_p$. Therefore, the centralizer of~$\Lie S$ is also
a (maximal) torus of~$\Lie G_{\rational_p}$. Because $\Lie S$ is
$\rational_p$-split, this implies that $\Lie G_{\rational_p}$ is
quasi-split.
 \end{proof}

\section{Which $\rational$-forms are
$\real$-universal?}

\begin{defn}
 Let us say that a Lie algebra~$\Lie G_\rational$
over~$\rational$ is \emph{universal for real representations}
(or simply \emph{$\real$-universal}, for short) if every real
representation of~$\Lie G_\rational$ has a $\rational$-form.
 \end{defn}

We have shown that every semisimple real Lie algebra has an
$\real$-universal $\rational$-form \see{EGoodQform}.
Furthermore, our construction yields an example that is
$\rational(i)$-split \fullsee{ConstructQform}{Qisplit}. In fact,
results of J.~Tits \cite{Tits-IrredRepns} show it is often the
case that \emph{every} $\rational(i)$-split $\rational$-form
of~$\Lie G$ is $\real$-universal.

\begin{prop} \label{TitsAllButSO}
 Let $\Lie G$ be a compact, simple Lie algebra over~$\real$.
If $\Lie G$ has a $\rational$-form that splits over some
quadratic extension of~$\rational$, but is not $\real$-universal,
then either
 $$ \mbox{$\Lie G \iso \su(n)$, for some even $n \ge 4$,} $$
 or
 $$ \mbox{$\Lie G \iso \so(n)$, for some $n \not\equiv 3$, $5$
$\pmod{8}$} .$$
 \end{prop}

In this section, we show that the converse is true \cf{badGQ}.

Using Tits' approach, it should be possible to give an explicit
list of the $\rational$-forms of $\su(n)$ and $\so(n)$ that are
not $\real$-universal (except, perhaps, those that are triality
forms of type $^3 D_4$ or $^6 D_4$). The author
intends to attack this project in a future paper.

\begin{proof}[{\bf Proof of Prop.~\ref{TitsAllButSO}}]
 Let $F$ be a quadratic extension of~$\rational$, and suppose
$\Lie G_\rational$ is a $\rational$-form of~$\Lie G$ that splits
over~$F$, but is not $\real$-universal. (We remark that $F$ must
be an imaginary extension, because $\Lie G_\real = \Lie G$ is
compact, not split.) There is an irreducible
$\rational$-representation~$V$ of~$\Lie G_\rational$, such that
$V$ is reducible over~$\real$ \see{Qform<>RIrred}. Thus, we may
write $V_\real = W \oplus X$ (with $W$ and $X$ nontrivial).
From Cor.~\fullref{Q(i)split}{2irreds}, we see that $W_\complex$
and~$X_\complex$ are irreducible. 

Let $\lambda$ be the highest weight of~$W_\complex$, and let
$w$~be the longest element of the Weyl group of~$\Lie G$. Because
$W_\complex$ is irreducible, and has an $\real$-form
(namely,~$W$), it must be the case that 
 \begin{itemize}
 \item $w(\lambda) = -\lambda$,
 and
 \item when~$\lambda$ is expressed as a linear combination of
simple roots, the sum of the coefficients is an
integer
 \end{itemize}
 (see \cite[Prop.~6.1 and following comments]{Tits-IrredRepns}).
Because $W_\complex$ has no $\rational$-form (any such
$\rational$-form would be isomorphic to a proper submodule
of~$V$), it must be the case that
 \begin{itemize}
 \item $\lambda$ is not an integral linear combination of
roots
 \end{itemize}
 (see \cite[Thm.~3.3]{Tits-IrredRepns}). Now Lem.~\ref{SUSO<>wt}
below yields the desired conclusion.
 \end{proof}

The following observation is obtained by inspection of a list
of the fundamental dominant weights of the complex simple Lie
algebras of each type $A_\ell, B_\ell, C_\ell, D_\ell, E_6, E_7,
E_8, F_4, G_2$. Such a list appears in \cite[Table~1,
p.~69]{HumphreysLieAlg}.

\begin{lem} \label{SUSO<>wt}
 Suppose $\Lie G$ is a simple Lie algebra
over~$\complex$, and let $w$ be the longest element of the Weyl
group of~$\Lie G$. There is a dominant weight~$\lambda$ of~$\Lie
G$, such that
 \begin{enumerate}
 \item $w(\lambda) = -\lambda$, and
 \item when~$\lambda$ is expressed as a linear combination of
simple roots, the sum of the coefficients is an
integer,
 and
 \item $\lambda$ is not an integral linear combination of roots
 \end{enumerate}
 if and only if either
 \begin{enumerate} \renewcommand{\theenumi}{\alph{enumi}}
 \item $\Lie G$ is of type~$A_\ell$, with $\ell$~odd
{\upshape(}and $\ell \ge 3${\upshape)}, or
 \item $\Lie G$ is of type~$B_\ell$, with $\ell \equiv 3$ or~$4$
$\pmod{4}$ {\upshape(}and $\ell \ge 3${\upshape)},
 or
 \item $\Lie G$ is of type~$D_\ell$ {\upshape(}and $\ell \ge
3${\upshape)}.
 \end{enumerate}
 \end{lem}

Combining Lem.~\ref{SUSO<>wt} with the following result
establishes the converse of Prop.~\ref{TitsAllButSO}.

\begin{prop} \label{badGQ}
 Suppose
 \begin{itemize}
 \item $G$ is a compact, real, semisimple Lie group,
 \item $\Lie G$ is the Lie algebra of~$G$;
 \item $\Lie T$ is a maximal torus of~$\Lie G$;
 \item $\Phi$~is the root system of~$\Lie G$ {\upshape(}with
respect to~$\Lie T${\upshape)}
 \item $w \in G$~is a representative of the longest element of
the Weyl group of~$G$;
 \item $V$~is an irreducible $\complex$-representation
of~$G$;
 and
 \item $\lambda$~is the highest weight of~$V$.
 \end{itemize}
 If
 \begin{itemize}
 \item $w(\lambda) = -\lambda$;
 \item when~$\lambda$ is expressed as a linear combination of
simple roots, the sum of the coefficients is an
integer;
 and
 \item $\lambda \notin \langle \Phi \rangle$ {\upshape(}that
is, $\lambda$ is not an integral linear combination of
roots{\upshape)};
 \end{itemize}
 then there exist
 \begin{enumerate}
 \item \label{VR}
 a real form~$V_\real$ of~$V$;
 and
 \item \label{GQ}
 a $\rational$-form~$\Lie G'_\rational$ of~$\Lie G$;
 \end{enumerate}
 such that
 \begin{enumerate} \renewcommand{\theenumi}{\alph{enumi}}
 \item \label{GQsplit}
 $\Lie G'_\rational$ splits over~$\Qi$;
 and
 \item \label{noVQ}
 $V_\real$ does not have a $\rational$-form
{\upshape(}with respect to~$\Lie G'_\rational${\upshape)}.
 \end{enumerate}
 \end{prop}

\begin{proof}
 \pref{VR} Since $w(\lambda) = -\lambda$, we see that the
lowest weight of~$V$ is~$-\lambda$. Since the dual of~$V$ is
the irreducible $\Lie G$-module whose lowest weight
is~$-\lambda$, we conclude that $V$~is self-dual. That is,
$V$~is isomorphic to its dual. 

Because $\Lie G$ is compact, there is a $\Lie G$-invariant
Hermitian form on~$V$, so $V$~is conjugate-isomorphic to its
dual. Combining this with the conclusion of preceding
paragraph, we conclude that $V$~is conjugate-isomorphic to
itself. That is, $V$ is isomorphic to its conjugate.
Therefore $V \otimesR \complex \iso V \oplus V$ is the direct
sum of two isomorphic irreducibles.

Because the sum of the coefficients of~$\lambda$ is an integer,
$V$~has a real form~$V_\real$ \see{RealVsQuat}.
(We remark that any two real forms of~$V$ are isomorphic, so it
does not matter which one is chosen.)

\medskip

(\ref{GQ}\ref{GQsplit}) Let $\Delta$ be a base of~$\Phi$. The
difference of the highest weight and the lowest weight is a
sum of roots, so, because the highest weight is~$\lambda$ and
the lowest weight is~$-\lambda$, we may write
 $$ 2 \lambda = \sum_{\delta \in \Delta}^r a_\delta \delta, $$
 with each $a_\delta \in \integer$. Because $\lambda \notin
\langle \Phi \rangle$, there must be some $\tau \in \Delta$,
such that $a_{\tau}$~is odd.

Let 
 $$ \{h_\delta\}_{\delta \in \Delta} \cup
\{x_\alpha\}_{\alpha \in \Phi} $$
 be the usual Chevalley basis of~$\Lie G_\complex$, so that
 $$ \Lie G = \bigoplus_{\delta \in \Delta} i \real
h_\delta
 \oplus 
 \bigoplus_{\alpha \in \Phi}
 \bigl\{\, z x_\alpha + \conj{z} x_{-\alpha} \mid z \in
\complex \, \bigr\}
 .$$
 For each $\alpha \in \Phi$, we may write 
 $$\alpha = \sum_{\delta \in \Delta} c_\delta(\alpha) \,
\delta ,$$
 with each $c_\delta(\alpha) \in \integer$. In particular, we
have defined a function $c_\tau \colon \Phi \to \integer$.
Now let
 $$ \Lie G'_\Qi
 = \bigoplus_{\delta \in \Delta} \Qi h_\delta
 \oplus \bigoplus_{\alpha \in \Phi}
 (\sqrt{3})^{c_\tau(\alpha)} \,
 \bigl(\Qi x_\alpha + \Qi x_{-\alpha} \bigr)
 $$
 and
 \begin{align*}
 \Lie G'_\rational
 &= \Lie G'_\Qi \cap \Lie G \\
 &= \bigoplus_{\delta \in \Delta} i \rational h_\delta
 \oplus \bigoplus_{\alpha \in \Phi}
 (\sqrt{3})^{c_\tau(\alpha)} \,
 \bigl\{\, z x_\alpha + \conj{z} x_{-\alpha} \mid z \in \Qi
\, \bigr\}
 .
 \end{align*}

Because 
 $$\Qi (\sqrt{3})^{-r} = 3^r\Qi (\sqrt{3})^{-r} = \Qi
(\sqrt{3})^r ,$$
 we see that $\Lie G'_\Qi$ is the $\Qi$-span of the Chevalley
basis
 $$ \{h_\delta\}_{\delta \in \Delta} \cup 
 \left\{\, (\sqrt{3})^{c_\tau(\alpha)} \, x_\alpha \,
\right\}_{\alpha \in \Phi} ,$$
 so $\Lie G'_\Qi$ is a split $\Qi$-form of~$\Lie G_\complex$.
From this (and because it is obvious that the $\real$-span of
$\Lie G'_\rational$ is all of~$\Lie G$), it follows that $\Lie
G'_\rational$ is a $\rational$-form of~$\Lie G$, such that
$\Lie G'_\rational$ splits over~$\Qi$.

\medskip

(\ref{GQ}\ref{noVQ}) Suppose $V_\real$~has a $\Lie
G'_\rational$-invariant $\rational$-form~$V'_\rational$.
(This will lead to a contradiction.) Then there is a $\Lie
G'_\rational$-equivariant  conjugate-linear involution
$\sigma' \colon V'_\Qi \to V'_\Qi$.

Let $\Lie G_\rational$ be the standard $\rational$-form
of~$\Lie G$ (obtained by replacing $\sqrt{3}$ with~$1$ in the
above construction). We know (from Cor.~\ref{RrepsQform})
that $V_\real$~has a $\Lie G_\rational$-invariant
$\rational$-form~$V_\rational$. Let $\sigma \colon V_\Qi \to
V_\Qi$ be the corresponding $\Lie G_\rational$-equivariant
conjugate-linear involution.

Because all highest-weight vectors of~$V$ are scalar
multiples of each other, there is no harm in assuming that
$(V'_\Qi)^{\lambda} \cap V_\Qi^{\lambda} \neq 0$. (Simply
replace $V_\Qi$ with $k V_\Qi$, for some $k \in \complex$.)
Then, because these are one-dimensional vector spaces
over~$\Qi$, we must have 
 $$(V'_\Qi)^{\lambda} = V_\Qi^{\lambda} .$$

Then, by induction on the length of~$\alpha$, it is easy to
see that if $\alpha$~is any integral combination of elements
of~$\Delta$, with non-negative coefficients, then
 $$ (V'_\Qi)^{\lambda - \alpha} = (\sqrt{3})^{c_\tau(\alpha)}
V_\Qi^{\lambda - \alpha} .$$
 In particular, because $c_\tau(2\lambda)$ is odd (this was
how~$\tau$ was chosen), we know that
 $$ (V'_\Qi)^{-\lambda} = \sqrt{3} \, V_\Qi^{-\lambda} .$$
 Because $\lambda(\Lie H'_\rational) = \lambda(\Lie H_\rational)
\subset i \rational$, we have
 $$ \mbox{$\sigma(V_\Qi^{-\lambda}) = V_\Qi^{\lambda}$
 and
 $\sigma'\bigl( (V'_\Qi)^{\lambda} \bigr) =
(V'_\Qi)^{-\lambda}$.} $$

Let $f = \sigma' \circ \sigma \colon V \to V$, so $f$~is
$\complex$-linear and $\Lie G$-equivariant; hence, $f$~is a
scalar, say $f(v) = kv$. Therefore
 \begin{align*}
  k \cdot V_\Qi^{-\lambda}
 &= f \bigl( V_\Qi^{-\lambda} \bigr) 
 = \sigma' \bigl( \sigma (V_\Qi^{-\lambda}) \bigr) 
 = \sigma' \bigl( V_\Qi^{\lambda} \bigr) \\
 &= \sigma' \bigl( (V'_\Qi)^{\lambda} \bigr) 
 = (V'_\Qi)^{-\lambda} 
 = \sqrt{3} \, V_\Qi^{-\lambda}
 ,
 \end{align*}
 so $k = \sqrt{3} k'$, for some $k' \in \Qi$.

We have
 $$ (\sigma')^2 = (f\circ \sigma)^2 = (k\sigma)^2 = (\sqrt{3}
k' \sigma)^2
 = 3 k' \cjg{k'}
 \neq 1 $$
 (because $3$~is not a sum of two rational squares).
 This contradicts the fact that $\sigma'$~is an involution.
 \end{proof}

\begin{eg}
 The direct sum of $\real$-universal Lie algebras need not be
$\real$-universal. For example, let $\Lie G_1 = \so(5)$ and
$\Lie G_2 = \so(11)$. For $j = 1,2$, the Lie algebra $(\Lie
G_j)_\complex$ has a (fundamental) dominant weight $\lambda_j$,
such that, when~$\lambda_j$ is expressed as a linear combination
of simple roots, the sum of the coefficients is a
half-integer. Then the weight $\lambda_1 \otimes \lambda_2$ of
$\Lie G_1 \oplus \Lie G_2$ satisfies the hypotheses of
Prop.~\ref{badGQ}, so $\Lie G_1 \oplus \Lie G_2$ has a
$\rational$-form that is not $\real$-universal. However, any
$\rational$-form of $\Lie G_1 \oplus \Lie G_2$ must be the
direct sum of $\rational$-forms of the factors and those, by
Prop.~\ref{TitsAllButSO}, are $\real$-universal.
 \end{eg}

The following useful observations are well known.
 $(\ref{RealVsQuat-irred} \Leftrightarrow
\ref{RealVsQuat-real})$ follows from Schur's Lemma.
 $(\ref{RealVsQuat-real} \Leftrightarrow
\ref{RealVsQuat-w2})$ follows from \cite[Lem.~79(b),
p.~226]{Steinberg}. 
 $(\ref{RealVsQuat-real} \Leftrightarrow
\ref{RealVsQuat-int})$ follows from  
\cite[Prop.~6.1]{Tits-IrredRepns}.
 $(\ref{RealVsQuat-irred} \Leftrightarrow
\ref{RealVsQuat-w2})$ is implicit in the proof on pp.~137--138
of \cite{Raghunathan-Qform}.

\begin{lem} \label{RealVsQuat}
 Suppose 
 \begin{itemize}
 \item $G$ is a compact, semisimple, real Lie group, 
 \item $V$ is an irreducible, self-dual, real representation
of~$G$, 
 \item $\lambda$ is the highest weight of~$V$,
 and 
 \item $w \in G$ is a representative of the longest element of
the Weyl group of~$G$.
 \end{itemize}
 Then:
 \begin{enumerate}
 \item We have $\End_G(V) \iso \real$ or~$\quaternion$.
 \item We have $w(\lambda) = -\lambda$ and $\lambda(w^2) = \pm
1$.
 \item We may write $\lambda$ as a linear combination of the
fundamental dominent weights of~$G$, and the sum~$s$ of the
coefficients in this linear combination is either an integer or
a half-integer.
 \item The following are equivalent:
 \begin{enumerate}
 \item \label{RealVsQuat-irred}
 $V \otimesR \complex$ is irreducible;
 \item \label{RealVsQuat-real}
 $ \End_G(V) \iso \real$;
 \item \label{RealVsQuat-w2}
 $\lambda(w^2) = 1$;
 \item \label{RealVsQuat-int}
 the sum $s$ is an integer.
 \end{enumerate}
 \end{enumerate}
 \end{lem}

\end{document}